# Evaluation of mechanical and energy properties for the phase field modeling of failure

Yuanfeng Yu[1], Xiaoya Zheng[1,*], Peng Li[2], Zhongzhou Zhang[2], Jinyou Xiao[1]

[1]School of Astronautics, Northwestern Polytechnical University, Xi'an 710072, China

[2]Xi'an Modern Control Technology Research Institute, Xi'an 710065, China

**Abstract**

In recent years, various phase field models have been developed in variational methods to simulate the failure of brittle solids. However, there is a lack of objective evaluation of the existing results, and in particular, there are few studies on model nonhomogeneous resolution, stress-strain linear elastic properties, and failure stress estimation. To compensate for the above gaps, the commonly used variational phase field model is systematically analyzed to solve the problem of evaluating the mechanics and energy properties of the model in this paper. The unified expression of the analytical solution and the nonhomogeneous solution under specific boundary conditions is analyzed and verified. Additionally, we theoretically analyze the energy properties of the phase field model and study the influence of the critical strain energy on the damage field, stress and strain of different models. Finally, the effect of different effective material parameters on the material failure stress is analyzed, the stress estimation formula for the failure of mode I of various phase field models is given, and the accuracy of the theory is verified by some examples.

**Keywords:** Phase field method, Variational method, Brittle fracture, Crack propagation, Critical strain energy

---

[*]**Corresponding author**

E-mail: zhengxy_8042@nwpu.edu.cn

## Nomenclature

| | |
|---|---|
| $G_c$ | critical energy release rate |
| $\gamma$ | surface density function |
| $\Psi$ | strain energy function |
| $\lambda, \mu$ | lame constant |
| $f$ | energy driving force |
| $\mathcal{H}_c$ | maximum historical reference energy |
| $d_{\text{hom}}$ | homogenous phase field solution |
| $d_{\text{nonhom}}$ | nonhomogenous phase field solution |
| $\sigma_e$ | elastic limit stress |
| $\sigma_c$ | peak stress |
| $d_c$ | critical phase field value |
| $\mathcal{H}_c$ | critical strain energy |
| $G_c^{\text{eff}}$ | effective critical energy release rate |
| $a_0$ | initial crack length |
| $a_0^{\text{eff}}$ | effective initial crack length |
| $\sigma_y^{0,b}$ | maximum tensile stress |
| $\sigma_y^{0,c}$ | critical failure stress |
| $E_0$ | elastic modulus |
| $v$ | Poisson's ratio |
| $\eta$ | critical stress scale factor |
| CDM | continuous damage mechanics |
| FM | fracture mechanics |
| VCCT | virtual crack closure technique |

| XFEM | extended finite element method |
| PD | peridynamics |
| PFM | phase field method |
| FEM | finite element method |
| IEM | interface element method |

# 1 Introduction

The failure of the material is a complex process of gradual change, and the stress and strain of the material continue to accumulate under the action of the load. After reaching the failure point of the material, damage begins to occur, and microcracks are formed inside the material. With the continuous increase in load, the cracks continue to expand and gradually form macrocracks, which eventually lead to material failure [1,2]. The analysis methods of the material failure process can be divided into continuous damage mechanics (CDM) [3] and fracture mechanics (FM) [4]. Continuous damage mechanics mainly studies the local deterioration of the macromechanical properties of materials caused by the initiation, convergence and expansion of microstructural defects in materials and reflects the influence of damage by decaying the constitutive relation of materials [5,6]. According to the different damage mechanisms, the corresponding damage criterion and damage evolution equation can be established. When the damage criterion is introduced into the material constitutive relation, the damage degree of the whole analysis process can be described. Fracture mechanics mainly studies the effect of macrocracks on material strength. This method simulates crack initiation and propagation to describe the accumulation of defects and failure [7].

The failure process of materials involves nonlinear stress and strain changes and crack initiation and propagation processes, which are difficult to analyze. Therefore, in the field of structural and mechanical research, numerical simulation has become an important means to analyze the failure process of materials. Numerical analysis methods can be divided into discrete methods and dispersion methods. Discrete methods include the virtual crack closure technique (VCCT) [8] and extended finite element method (XFEM) [9]. When VCCT is used to simulate crack propagation, it is necessary to arrange the crack in advance, and the mesh near the crack needs to be subdivided or redivided continuously. XFEM describes the discontinuous term by adding the enrichment term of displacement, but its convergence is poor. In addition, it is necessary to know the possible failure mode in advance and define the corresponding cracking direction after damage. Numerical methods of dispersion, such as peridynamics (PD) [10] and the phase field method (PFM) [11], disperse cracks and describe the crack propagation process. In the peridynamics, the deformation compatibility equation and force balance equation in solid mechanics are rewritten as integral equations to describe crack initiation, propagation and bifurcation [12,13]. However, it is difficult for this method to directly correspond to the constitutive model of the material, and it is also incompatible with the finite element method (FEM) in the numerical solution.

Research on the phase field method originated in the late 1990s and has become an important research direction since then [14]. In this method, a field variable (order parameter)—the phase field—is introduced to describe the transition between different physical fields in a system. The variable varies between 0 and 1, representing the intact state and complete failure state of the material, respectively, thereby describing the discontinuous nature of the crack. This method can be used to simulate complex fracture processes, such as crack initiation, propagation, bifurcation and intersection [15].

At present, the phase field fracture models for analyzing material failure are mainly divided into two types: the first type is based on the Ginzburg-Landau theory in physics [16,17]. The Ginzburg-Landau type evolution equation and nonconvex attenuation function of the unconstrained crack phase field are used to describe the damage evolution process of materials and solve the Ginzburg-Landau equation

$$\frac{\partial d}{\partial t} = -M_d \frac{\partial \psi}{\partial d} \quad (1)$$

The displacement field and damage field are obtained. In Eq. (1), $d$ is the order parameter, $\psi$ is the total free energy of the system with respect to the order parameter, and $M_d$ is the model parameter of the corresponding system. It can be seen from the Eq. (1) that the model is completely viscous in nature and is mainly suitable for dynamic fracture [11]. The second kind of variational method originating from fracture mechanics is an extension of the classical Griffith fracture theory [18]. The theory was proposed by Francfort and Marigo [19] and Bourdin et al. [20]. The variational method and convergence criterion are used to study the crack evolution process, which makes it easier for researchers to understand its mechanical concept and put it into engineering application. In this model, elastic potential energy and crack surface energy are used to describe the total potential energy of the system

$$\mathcal{E}(d) = \int_\Omega W(\nabla d) d\Omega + \mathcal{H}^1(S(d)) \quad (2)$$

By solving the variational problem of Eq. (2) about the displacement field $u$ and the phase field $d$, the force balance equation and damage evolution equation are obtained, and then by solving the two equations, the displacement field $u$ and the damage field $d$ are obtained, where $W(\nabla d)$ is the elastic energy density, and $\mathcal{H}^1(S(d))$ is the crack surface energy of the material fracture surface. The theory is further developed by Miehe et al. [11]. The unification of the two types of phase field models is established

$$\Pi_e(\nabla_s \boldsymbol{u}, d, \nabla d) = g(d)\psi_0(\nabla_s \boldsymbol{u}, d) + G_c \gamma(d, \nabla d) \tag{3}$$

$$\begin{aligned}
&\dot{d} - \frac{1}{\eta}\left\langle -G_c \delta_d \gamma - g'(d)\psi_0 \right\rangle_+ = 0 \\
&\Rightarrow \eta \dot{d} = \left\langle -G_c \delta_d \gamma - g'(d)\psi_0 \right\rangle_+ \\
&\Rightarrow \eta \dot{d} = \left\langle -\delta_d \Pi_e(\nabla_s \boldsymbol{u}, d, \nabla d) \right\rangle_+
\end{aligned} \tag{4}$$

where $g(d)$ is the degradation function, $\psi_0$ is the undamaged elastic energy density, $G_c$ is the critical energy release rate of the material, $\gamma(d, \nabla d)$ is the crack surface density function, $\eta$ is the viscosity coefficient, which is equivalent to the model parameters in Eq. (1), and $\langle x \rangle_\pm = (x \pm |x|)/2$ is the Macaulay bracket operator. At the same time, Miehe et al. [21] also creatively proposed the concept of a historical field, which makes the evolution of the crack phase field meet irreversible constraints. This lays a very important theoretical foundation for the study of the crack growth process under cyclic loading and paves an important way for the popularization and application of phase field theory so that it can solve complex crack growth problems. Based on Miehe's research, the phase field method has been improved and applied to different research fields. These include brittle fracture [22,23], plastic fracture [24,25], dynamic fracture [26,27], cohesive fracture [28,29] and material fatigue [30,31]. At the same time, the scope of research materials is also expanding, such as rock materials [32,33], polymers [34,35], concrete materials [36,37], and composite materials [38,39]. Compared with the existing methods, such as the XFEM method and interface element method (IEM), the phase field method has the following advantages in the study of the material failure process, which enables it to describe discontinuous cracks and simulate complex crack propagation.

  (1) The phase field method is used to solve the crack propagation process based on the principle of energy minimization and can be used to automatically predict crack nucleation, propagation and combination without defining the direction of crack propagation in advance;

(2) The phase field method can be used to predict crack initiation and bifurcation processes without additional criteria, which provides a new idea to simulate complex crack paths;

(3) The phase field method overcomes the mesh dependence problem in the process of the finite element method solving material damage and fracture.

For the crack density function in Eq. (3), there are generally three types, commonly used as AT1 and AT2 models [40]. For the AT1 model [41], the stress-strain relationship has a linear elastic stage and has an elastic limit stress; that is, damage will occur after a period of loading. For the AT2 model [42], the stress-strain relationship is always nonlinear, and there is no elastic limit; that is, damage occurs immediately when the load is applied, and when the material reaches the limit stress, the damage field has increased to $d = 0.25$, which affects the accuracy of the calculation results to a certain extent. There is no specific and detailed analysis for why this phenomenon occurs when different models are selected. In this context, this paper will conduct a theoretical analysis of the mechanical and energy properties of different phase field damage models, reveal the causes of this damage phenomenon, and study the mechanical change properties of different damage models to provide theoretical guidance for the selection of phase field damage models.

After introduction, this paper is organized as follows: Section 2 introduces the phase field theory, Section 3 deduces the homogenous solution, nonhomogenous solution and limit stress of the phase field theory, and analyzes the elastic and energy characteristics of the damage model, as well as the mechanical characteristics of the phase field damage model. Section 4 introduces the finite element realization process of the phase field model, Section 5 verifies the theoretical analysis by some examples, and Section 6 summarizes the research results of this paper.

**2 Phase field fracture model**

The phase field method is summarized in this section, and the governing equation of the phase field model is derived, which is convenient for the theoretical analysis in the next section.

*2.1 Fracture energy functional*

*2.1.1 Crack dissipation functional*

For a one-dimensional bar that breaks in the axial position $x=0$, where $\Gamma$ represents the crack surface completely destroyed, as shown in Fig. 1. An auxiliary variable $d(x) \in [0,1]$ can be introduced to describe the sharp crack, $d=0$ represents the intact state of the bar, and $d=1$ represents the completely broken state of the bar.

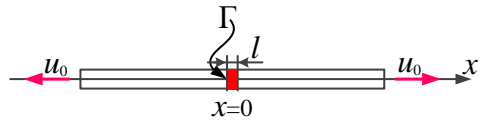

**Fig. 1.** 1D bar crack

Here, an exponential distribution function [11]

$$d(x) = e^{-|x|/l} \tag{5}$$

can be used to approximate the nonsmooth crack phase field to disperse the sharp cracks on the axial region $L$ of the bar. A parameter $l$ controls the dispersion width of the crack and affects the variation trend of the field variables, as shown in Fig. 2. When the parameter $l$ changes from large to small, the curve begins to change from "fat" to "thin". At that time $l \to 0$, the diffuse cracks recover to sharp cracks, corresponding to the red solid line in Fig. 2.

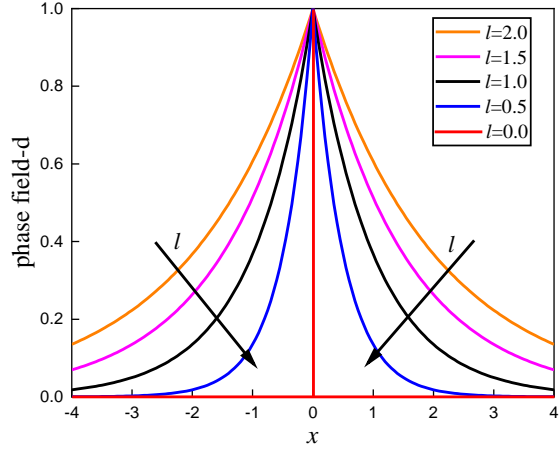

**Fig. 2.** Field variable $d$ varies with $l$

According to Eq. (5), one can further use functional

$$\Gamma(d) = \frac{1}{2l}\int_\Omega \{d^2 + l^2 d'^2\}dV \tag{6}$$

to characterize a sharp crack topology.

For the multidimensional case, the model area is $\Omega \subset R^\delta$, the spatial dimension is $\delta \in [1-3]$, and its surface is $\partial\Omega \subset R^{\delta-1}$, as shown in Fig. 3.

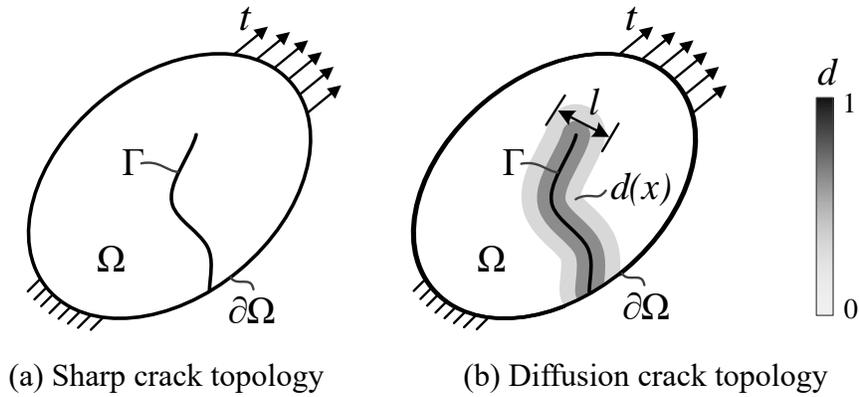

(a) Sharp crack topology  (b) Diffusion crack topology

**Fig. 3.** 2D crack topology

Here, the crack functional Eq. (6) is extended to the multidimensional case

$$\Gamma(d) = \int_\Omega \gamma(d, \nabla d)dV \tag{7}$$

The surface density function $\gamma$ and derivative $\delta_d\gamma$ of the crack per unit area (or length) are

$$\gamma(d,\nabla d)=\frac{1}{2l}d^2+\frac{l}{2}|\nabla d|^2, \quad \delta_d\gamma=\frac{1}{l}d-l\Delta d.  \tag{8}$$

For the crack area density function in Eq. (8), the crack density function can be expressed in the following general form:

$$\gamma(d,\nabla d)=\frac{1}{4c_w}\left(\frac{1}{l}w(d)+l|\nabla d|^2\right), \quad \delta_d\gamma=\frac{1}{4c_w}\left(\frac{1}{l}w'(d)-2l\Delta d\right) \tag{9}$$

Here, the parameters $c_w=\int_0^1\sqrt{w(\hat{d})}d\hat{d}$ and $w(d)$ determine the distribution of the crack phase field [43],

$$w(d)=\xi d+(1-\xi)d^2 \in[0,1] \quad \forall d \in[0,1] \tag{10}$$

In Eq. (10), the parameters $\xi \in [0,2]$ and geometric crack function are shown in Table 1.

| $w(d)$ | $\xi$ | $c_w$ | $d(x)$ |
|---|---|---|---|
| $d^2$ | 0 | $\frac{1}{2}$ | $\exp\left(-\frac{|x|}{l}\right)$ |
| $d$ | 1 | $\frac{2}{3}$ | $\left(1-\frac{|x|}{2l}\right)^2$ |
| $2d-d^2$ | 2 | $\frac{\pi}{4}$ | $1-\sin\left(\frac{|x|}{l}\right)$ |

Table 1. Crack function $w(d)$ and form of geometric crack function $d(x)$

For the phase field fracture problem in Fig. 3, according to the fracture surface function $\Gamma(d)$ introduced in Eq. (7), the work required to generate a new crack topology can be defined

$$\Pi_c(d) = \int_\Omega G_c \gamma(d, \nabla d) dV \tag{11}$$

Here, $\gamma(d, \nabla d)$ is the crack surface density function defined in Eq. (9), and $G_c$ is the Griffith-type critical energy release rate, which is the work required to produce a new unit length crack [18].

*2.1.2 Strain energy functional*

According to linear elastic theory, the global energy storage functional can be expressed as

$$\Pi_e(\boldsymbol{u}, d) = \int_B \psi(\varepsilon(\boldsymbol{u}), d) dV, \tag{12}$$

The energy storage function $\psi(\varepsilon, d)$ describes the volume energy stored in the unit volume of the solid, and the energy reduction caused by fracture can be expressed as

$$\psi(\varepsilon, d) = [g(d) + k] \psi_0(\varepsilon). \tag{13}$$

In the above formula, $g(d)$ is the degradation function, which satisfies the property

$$\begin{cases} g(0) = 1, \ g(1) = 0. \\ g(d) > 0, \ for \ d \neq 1 \\ g'(d) < 0, \ for \ d \neq 1 \\ g'(0) < 0, \ g'(1) = 0. \end{cases} \tag{14}$$

The parameter $k$ is the residual stiffness coefficient to ensure numerical stability. $\psi_0$ is the strain energy function of the undamaged material. For isotropic material, it can be expressed as

$$\psi_0(\varepsilon) = \lambda \text{tr}^2[\varepsilon]/2 + \mu \text{tr}[\varepsilon^2]. \tag{15}$$

with lame constant $\lambda > 0$ and $\mu > 0$. To avoid compression damage, the strain energy needs to be decomposed into the positive part $\psi_0^+(\varepsilon)$ and the negative part $\psi_0^-(\varepsilon)$,

$$\psi_0(\varepsilon) = \psi_0^+(\varepsilon) + \psi_0^-(\varepsilon).  \tag{16}$$

According to the spectral decomposition method of Miehe et al. [21], $\psi_0^+(\varepsilon)$ and $\psi_0^-(\varepsilon)$ can be expressed as

$$\psi_0^+(\varepsilon) = \lambda \langle \text{tr}(\varepsilon) \rangle_+^2 / 2 + \mu \text{tr}(\varepsilon_+^2),\ \psi_0^-(\varepsilon) = \lambda \langle \text{tr}(\varepsilon) \rangle_-^2 / 2 + \mu \text{tr}(\varepsilon_-^2). \tag{17}$$

Therefore, the stress field can be calculated from Eq. (13)

$$\boldsymbol{\sigma} = \partial_\varepsilon \psi(\varepsilon, d) = (g(d) + k)\left[ \lambda \langle \text{tr}(\varepsilon) \rangle_+ \boldsymbol{I} + 2\mu \varepsilon_+ \right] + \left[ \lambda \langle \text{tr}(\varepsilon) \rangle_- \boldsymbol{I} + 2\mu \varepsilon_- \right]. \tag{18}$$

## 2.2 Variational phase field model

According to the crack dissipation functional Eq. (11) and energy storage functional Eq. (12), the total energy functional is introduced

$$\Pi(\boldsymbol{u}, d) = \int_{\mathcal{B}} \psi(\varepsilon(\boldsymbol{u}), d) dV + \int_{\Omega} G_c \gamma(d, \nabla d) dV - \int_{\mathcal{B}} \boldsymbol{b} \cdot \boldsymbol{u} dV - \int_{\partial \mathcal{B}} \boldsymbol{t} \cdot \boldsymbol{u} dA. \tag{19}$$

where $b$ is the volume force in area $\mathcal{B}$ and $t$ is the surface external force on surface $\partial \mathcal{B}$.

The displacement field and phase field $(\boldsymbol{u}, d)$ can be determined by the following minimization problem.

$$(\boldsymbol{u}(x), d(x)) = \text{Arg}\{\inf \Pi(\boldsymbol{u}, d)\}\ \dot{d} \geq 0, d \in [0, 1] \tag{20}$$

According to the boundary condition of the total energy functional

$$\begin{cases} \delta \Pi = 0,\ \delta d > 0 \\ \delta \Pi > 0,\ \delta d = 0 \end{cases} \tag{21}$$

one can obtain the following strong form of the force balance equation:

$$\begin{aligned} \nabla \cdot \boldsymbol{\sigma} + b &= 0,\ \ \text{in } \mathcal{B} \\ \boldsymbol{\sigma} \cdot \boldsymbol{n} &= t,\ \ \ \ \text{on } \partial \mathcal{B} \end{aligned} \tag{22}$$

and phase field evolution equation

$$\begin{cases} f - G_c \delta_d \gamma = 0 \ \dot{d} > 0 \\ f - G_c \delta_d \gamma < 0 \ \dot{d} = 0 \end{cases}, \quad in\ \Omega$$

$$\frac{\partial \gamma}{\partial \nabla d} \cdot n_\Omega = 0, \qquad on\ \partial\Omega$$

(23)

The first condition in Eq. (23) is an equation when damage occurs $(\dot{d} > 0)$, and when it is stable $(\dot{d} = 0)$, this equation is an inequality equation. Here, the energy driving force $f$ conjugated to the energy functional $\psi(\varepsilon, d)$ is

$$f := -\frac{\partial \psi}{\partial d} = -g'(d)\psi_0^+ \tag{24}$$

To ensure the irreversibility of damage evolution, a positive maximum historical reference energy is introduced

$$\mathcal{H}(x,t) := \max_{s \in [0,t]}\left[\psi_0^+(\varepsilon(x,t)), \mathcal{H}_0\right],$$

$$\mathcal{H}_0 = -\frac{w'(0)}{4g'(0)c_w l} G_c. \tag{25}$$

Here, substituting Eq. (25) for $\psi_0^+$ in Eq. (24), one gets

$$f := -g'(d)\mathcal{H}. \tag{26}$$

**3 Analysis of mechanical and energy properties of the model**

*3.1 Research on analytical solution*

To illustrate some properties of the phase field fracture model, the analytical solution of the boundary value problem for isotropic materials under the small deformation theory is studied in this section. Here, one mainly analyzes the phase field fracture problem in the one-dimensional region $x \in [-L, L]$, and its properties can also be extended to multidimensional problems. It is assumed that there is no compressive strain field in this region, that is, $\psi_0^- = 0, \psi_0 = \frac{1}{2}E_0\varepsilon^2$.

Based on the above assumptions, the one-dimensional problem can be expressed as

$$\begin{cases} \dfrac{d\sigma}{dx} = 0 \\ \dfrac{G_c}{4c_w}\left[\dfrac{1}{l}w'(d) - 2l\Delta d\right] + g'(d)\psi_0 = 0 \end{cases} \tag{27}$$

where the stress is $\sigma = g(d)E_0\varepsilon$.

To simplify the mathematical expression, a monotonically increasing function $h(d)$ is introduced

$$h(d) = \frac{1}{g(d)} - 1 \tag{28}$$

As a result, it has been achieved

$$g(d) := \frac{1}{1+h(d)}, \quad g'(d) = -g^2(d)h'(d) < 0 \tag{29}$$

In Eq. (29), $h'(d) > 0$, Eq. (27) becomes

$$\sigma^2 h'(d) - A\left[w'(d) - 2l^2\Delta d\right] = 0 \tag{30}$$

Where parameter $A$ is $A = \dfrac{E_0 G_c}{2c_w l}$.

*3.1.1 Homogenous solution*

By neglecting the spatial gradient $\Delta d$ of the phase field $d$, a homogenous solution of the one-dimensional problem can be obtained. At this time, Eq. (30) becomes

$$\sigma^2 h'(d) - Aw'(d) = 0. \tag{31}$$

According to Eq. (31), the homogenous solution $d_{\text{hom}}(\sigma)$ can be obtained, and the stress field and strain field can also be obtained.

$$\sigma = \sqrt{A\frac{w'(d)}{h'(d)}}, \quad \varepsilon = \frac{1}{E_0}\sqrt{-A\frac{w'(d)}{g'(d)}}. \tag{32}$$

Here, the critical phase field values of different models can be obtained by calculating the extremum in the formula. For example, when the attenuation function is taken, the critical phase field values can be obtained

Here, by calculating the extreme value of the stress in formula (32), the critical phase field value ($d_c$) of different models can be obtained. For example, when the degradation function is taken as $g(d) = (1-d)^2$, one can obtain

$$d_c = \begin{cases} 0.25, & \xi = 0 \\ 0, & \xi = 1 \\ 0, & \xi = 2 \end{cases} \tag{33}$$

*3.1.2 Nonhomogenous solution*

The homogenous solution of the above analysis is only stable when the state is $\varepsilon \leq \varepsilon_c$ [44], but when the state is $\varepsilon > \varepsilon_c$, damage evolution will occur, leading to crack propagation, and the crack region will be controlled within the range $[-D, D]$ by the length parameter $l$. At this time, according to the conditions

$$\begin{cases} h(x = \pm D) = h(d = 0) = 0 \\ d(x = \pm D) = 0 \\ \nabla d(x = \pm D) = 0 \end{cases} \tag{34}$$

the general form of a nonhomogenous solution can be derived.

After integrating Eq. (30), It becomes

$$\sigma^2 h(d) - A\left[w(d) - 2l^2 \nabla d\right] = 0 \tag{35}$$

In the region far away from the crack, the phase field gradient disappears. At this time, the homogenous phase field value can be calculated according to the boundary conditions, namely, $\lim_{x \to \pm\infty} d(x) = d_{hom}(\sigma)$.

Multiplying Eq. (30) by $\dfrac{d\hat{d}}{dx}$ and integrating it, one gets

$$\frac{d}{dx}\left\{\sigma^2 h(\hat{d}) - A\left[w(\hat{d}) - l^2\left(\frac{d\hat{d}}{dx}\right)^2\right]\right\} = 0 \tag{36}$$

Therefore, it can be deduced that

$$\frac{d\hat{d}}{dx} = -\operatorname{sgn}(x)\sqrt{\frac{1}{Al^2}\left[Aw(\hat{d}) - \sigma^2 h(\hat{d}) - a\right]} \tag{37}$$

Where $a$ is a constant. For example, when the parameters $w(\hat{d})$ and $g(\hat{d})$ are

$$\begin{cases} w(\hat{d}) = \hat{d}^2 \\ g(\hat{d}) = (1-\hat{d})^2 \end{cases} \tag{38}$$

A value of the constant $a$ can be calculated as

$$a = -\frac{\sigma^2}{(1-d_{hom})^2} + \frac{E_0 G_c}{l} d_{hom}^2 + \sigma^2, \quad d_{hom} = \frac{E_0 l \varepsilon^2}{E_0 l \varepsilon^2 + G_c} \tag{39}$$

For the case of complete fracture, it can be obtained from Eq. (35)

$$\begin{cases} \sigma = 0, \\ d_{hom}(0) = 0. \end{cases} \tag{40}$$

As a result, Eq. (37) becomes

$$\frac{d\hat{d}}{dx} = -\operatorname{sgn}(x)\frac{\hat{d}}{l} \tag{41}$$

Assuming that a fracture occurs at $x=0$, that is, $d(0)=1$, it can be obtained from Eq. (41)

$$\hat{d}(x) = d_{nonhom} = e^{-|x|/l} \tag{42}$$

It can be found that Eq. (42) and (5), which indicates that the nonhomogeneous solution $d_{nonhom}$ that satisfies the specific boundary conditions is the geometric crack function of the phase field model.

Similarly, at that time $w(\hat{d}) = \hat{d}, 2\hat{d} - \hat{d}^2$, it can also be deduced

$$\hat{d}(x) = d_{nonhom} = \begin{cases} \left(1 - \frac{|x|}{2l}\right)^2, & w(\hat{d}) = \hat{d} \\ 1 - \sin\left(\frac{|x|}{l}\right), & w(\hat{d}) = 2\hat{d} - \hat{d}^2 \end{cases} \tag{43}$$

Comparing the nonhomogeneous solutions (42)-(43) that meet the special boundary conditions with the geometric crack function $d(x)$ in Table 1, it can be found that they are the same, which verifies the accuracy of the nonhomogeneous solution.

*3.1.3 Stress properties analysis*

1) Quadratic crack function

Taking $\xi = 0$, $w(d) = d$, and $g(d)$ taking the quadratic degradation function, according to the properties of Eq. (33), $d_c = 0.25$, one can obtain

$$\sigma_c = \frac{9}{16}\sqrt{\frac{E_0 G_c}{3l}}, \varepsilon_c = \sqrt{\frac{G_c}{3E_0 l}} \tag{44}$$

Further, the stress value can be obtained from Eq. (32)

$$\sigma = \left(\frac{G_c}{G_c + E_0 l \varepsilon^2}\right)^2 E_0 \varepsilon, \, d = \frac{E_0 l \varepsilon^2}{G_c + E_0 l \varepsilon^2} \tag{45}$$

It is found that there will be no elastic limit value for the stress function.

2) Linear crack function

Similarly, taking $\xi = 1$, $w(d) = d$, and $g(d)$ taking the quadratic degradation function, according to the properties of Eq. (33), $d_c = 0$, one can obtain

$$\sigma_c = \sqrt{\frac{3 E_0 G_c}{8l}}, \, \varepsilon_c = \sqrt{\frac{3 G_c}{8 E_0 l}} \tag{46}$$

and the following properties:

(a) Linear elastic stage: $\varepsilon \leq \bar{\varepsilon}$

$$\sigma = E_0 \varepsilon, \, d = 0 \tag{47}$$

Here, $\bar{\varepsilon} = \sqrt{\dfrac{3 G_c}{8 E_0 l}}$. at this time, the elastic limit stress is

$$\sigma_e = \sqrt{\frac{3 E_0 G_c}{8l}} \tag{48}$$

(b) Inelastic stage: $\varepsilon > \bar{\varepsilon}$

$$\sigma = \left(\frac{3 G_c}{8 E_0 l \varepsilon^2}\right)^2 E_0 \varepsilon, \, d = \frac{8 E_0 l \varepsilon^2 - 3 G_c}{8 E_0 l \varepsilon^2} \tag{49}$$

It is found that there will be an elastic limit value for the stress function.

(c) Mixed crack function

Similarly, taking $\xi = 2$, $w(d) = 2d - d^2$, and $g(d)$ taking the quadratic degradation function, according to the properties of Eq. (33), $d_c = 0$, one can obtain

$$\sigma_c = \sqrt{\frac{2E_0 G_c}{\pi l}}, \varepsilon_c = \sqrt{\frac{2G_c}{\pi E_0 l}} \tag{50}$$

and the following properties:

    (a)    Linear elastic stage: $\varepsilon \leq \bar{\varepsilon}$

$$\sigma = E_0 \varepsilon, \ d = 0 \tag{51}$$

Here, $\bar{\varepsilon} = \sqrt{\frac{2G_c}{\pi E_0 l}}$. at this time, the elastic limit stress is

$$\sigma_e = \sqrt{\frac{2E_0 G_c}{\pi l}} \tag{52}$$

    (b)    Inelastic stage: $\varepsilon > \bar{\varepsilon}$

$$\sigma = 0, \ d = 1 \tag{53}$$

It is found that there will also be an elastic limit value for the stress function.

According to the above calculation, the results of the stress and strain of different damage models are summarized, as shown in Table 2.

| Model | $w(d)$ | $g(d)$ | $d_c$ | $\sigma_e$ | $\sigma_c$ | $\varepsilon_c$ |
|---|---|---|---|---|---|---|
| 1 | $d^2$ | $(1-d)^2$ | 0.25 | 0 | $\frac{9}{16}\sqrt{\frac{E_0 G_c}{3l}}$ | $\sqrt{\frac{G_c}{3E_0 l}}$ |
| 2 | $d$ | $(1-d)^2$ | 0 | $\sqrt{\frac{3E_0 G_c}{8l}}$ | $\sqrt{\frac{3E_0 G_c}{8l}}$ | $\sqrt{\frac{3G_c}{8E_0 l}}$ |
| 3 | $2d-d^2$ | $(1-d)^2$ | 0 | $\sqrt{\frac{2E_0 G_c}{\pi l}}$ | $\sqrt{\frac{2E_0 G_c}{\pi l}}$ | $\sqrt{\frac{2G_c}{\pi E_0 l}}$ |

**Table 2.** Stress and strain properties

In the table: $\sigma_e$ is the elastic limit stress, $\sigma_c$ is the peak stress.

### *3.2 Analysis of the linear elastic and energy properties of the model*

It can be seen from Section 3.1 that in case $\xi = 0$, when the material begins to load, damage begins to occur, and there is no elastic limit stress in the model; when the stress increases to the peak stress, the material damage field increases to $d_c = 0.25$. However, in case $\xi = 1, 2$, the critical phase field value of the model is $d_c = 0$, which prevents the material from damaging at the beginning of loading so that the stress-strain relationship of the material maintains linear elastic changes at the initial stage, and the model has elastic limit stress. For different damage models, why does this phenomenon occur? Currently, there is no reasonable explanation for this problem. Therefore, this section will analyze the linear elastic properties and energy properties of the phase field model to clearly explain the reasons for this problem.

The phase field evolution equation is obtained from Section 2.2

$$\frac{G_c}{4c_w l}\left[\left(\xi + 2(1-\xi)d\right) - 2l^2 \Delta d\right] + g'(d)\mathcal{H} = 0 \tag{54}$$

Neglecting the effect of the spatial gradient of the phase field $d$, it can be obtained from Eq. (54)

$$\frac{G_c}{4c_w l}\left(\xi + 2(1-\xi)d\right) + g'(d)\mathcal{H} = 0 \tag{55}$$

and considering the initial damage condition $d = 0$, it can be further obtained

$$\frac{G_c}{4c_w l}\xi + g'(0)\mathcal{H}_c = 0 \tag{56}$$

Here, $\mathcal{H}_c$ is defined as the critical strain energy at the beginning of damage. From Eq. (56), the critical strain energy can be expressed as

$$\mathcal{H}_c = -\frac{\xi}{4g'(0)c_w l}G_c \propto G_c \tag{57}$$

It can be found from Eq. (57) that the critical strain energy $\mathcal{H}_c$ is directly proportional to the critical energy release rate $G_c$. When the phase field model is fixed, the greater the critical energy release rate $G_c$ is, the greater the critical strain energy $\mathcal{H}_c$, and vice versa.

To analyze the energy properties of the phase field fracture model, it is assumed that there is no compressive strain field in the region, namely, $\psi_0^- = 0, \psi_0 = \frac{1}{2}E_0\varepsilon^2$. In case $\xi = 0$, the critical strain energy can be obtained from Eq. (57)

$$\mathcal{H}_c = 0 \tag{58}$$

and elastic limit stress

$$\sigma_e = 0 \tag{59}$$

Similarly, in case $\xi = 1, 2$, the critical strain energy can be obtained

$$\mathcal{H}_c = \begin{cases} -\dfrac{3G_c}{8g'(0)l} > 0, & \xi = 1 \\ -\dfrac{2G_c}{\pi g'(0)l} > 0, & \xi = 2 \end{cases} \tag{60}$$

and elastic limit stress

$$\sigma_e = \begin{cases} \sqrt{-\dfrac{3G_c E_0}{4g'(0)l}}, & \xi = 1 \\ \sqrt{-\dfrac{4G_c E_0}{\pi g'(0)l}}, & \xi = 2 \end{cases} \tag{61}$$

Here, taking the quadratic degradation function $g(d) = (1-d)^2$, the damage field is obtained from Eq. (55)

$$d = \frac{\langle \mathcal{H} - \mathcal{H}_c \rangle_+}{\mathcal{H} + \dfrac{(1-\xi)G_c}{4c_w l}}, \quad \mathcal{H}_c = \frac{\xi G_c}{8c_w l} \tag{62}$$

Therefore, one can get

$$d = \begin{cases} 0, & \mathcal{H} \leq \mathcal{H}_c \\ \dfrac{\mathcal{H} - \mathcal{H}_c}{\mathcal{H} + \dfrac{(1-\xi)G_c}{4c_w l}}, & \mathcal{H} > \mathcal{H}_c \end{cases} \tag{63}$$

And stress properties

$$g(d) = \begin{cases} 1, & \mathcal{H} \leq \mathcal{H}_c \\ g(d), & \mathcal{H} > \mathcal{H}_c \end{cases}$$

$$\Rightarrow \sigma = g(d) E_0 \varepsilon = \begin{cases} E_0 \varepsilon, & \mathcal{H} \leq \mathcal{H}_c \\ g(d) E_0 \varepsilon, & \mathcal{H} > \mathcal{H}_c \end{cases} \tag{64}$$

It can be seen from Eq. (63) and (64) When the elastic strain energy of the material is less than or equal to the critical strain energy $\mathcal{H}_c$, no damage will occur, and the stress-strain relationship of the material will maintain a linear elastic change. Only when the elastic strain energy of the material is greater than the critical strain energy $\mathcal{H}_c$ will the material start to damage, and then crack propagation will occur, which shows that the critical strain energy $\mathcal{H}_c$ affects the time when the material starts to be damaged and then affects the stress-strain properties of the material.

Comparing Eq. (58) and (60), it can be seen that only when the critical strain energy $\mathcal{H}_c$ satisfies the relationship

$$\mathcal{H}_c > 0 \tag{65}$$

the elastic limit stress $\sigma_e$ exists in the model

$$\sigma_e = \sqrt{-\dfrac{\xi G_c E_0}{2 g'(0) c_w l}} \propto \sqrt{\dfrac{G_c E_0}{l}}. \tag{66}$$

Because there is no critical strain energy in the model in case $\xi = 0$, the damage model does not have an elastic limit stress, but in case $\xi = 1, 2$, the model has a critical strain energy, resulting in an elastic limit stress in the model.

To further verify the influence of the critical strain energy on the damage field and stress-strain field, an additional critical strain energy is introduced into the model $\xi = 0$ [45,46]

$$\mathcal{H}_c = \frac{G_c}{2l} \tag{67}$$

Then, the damage field and elastic limit stress can be deduced

$$d = \frac{\langle \mathcal{H} - \mathcal{H}_c \rangle_+}{\langle \mathcal{H} - \mathcal{H}_c \rangle_+ + \frac{G_c}{2l}} \tag{68}$$

$$\sigma_e = E_0 \varepsilon_c = \sqrt{\frac{E_0 G_c}{l}}$$

It can be found from the formula that at that time, the material will not be damaged, and the stress-strain relationship remains linear elastic. Only at that time, the material will be damaged. The above analysis further verifies that the critical strain energy affects the time when the material begins to damage and only when the critical strain energy exists in the model; that is, the elastic limit stress exists in the damage model.

Incorporating Eq. (68)₁ into Eq. (64), it can be found that when the strain field meets the condition $\mathcal{H} \leq \mathcal{H}_c$, the material will not be damaged, and the stress-strain relationship maintains a linear elastic relationship. Only when the strain field meets the condition $\mathcal{H} > \mathcal{H}_c$ will the material be damaged. The above analysis further verifies that the critical strain energy affects the time when the material starts to be damaged, and only when the critical strain energy $\mathcal{H}_c$ exists in the model, that is, $\mathcal{H}_c > 0$, will the model have the elastic limit stress.

### *3.3 Analysis of mechanical properties of the model*

When there is a crack in the material, due to the influence of the spatial gradient of the phase field $d$, the analytical solutions of the damage field, stress and critical strain cannot be obtained according to Eq. (55). At the same time, the effective length and effective energy release rate of different models also affect the mechanical properties of the model, so the results of stress and strain during material failure will be different from those in Section 3.1. Therefore, this section will study the mechanical properties of different models in the presence of cracks.

Under the influence of mesh size, the effective critical energy release rate $G_c^{eff}$ is defined as [47]

$$G_c^{eff} = G_c\left(1+\frac{h}{4c_w l}\right) \tag{69}$$

Here, $h$ is the minimum element feature length. It can be seen from Eq. (69) that when the parameters $h$, $l$ and $G_c$ are the same, $G_c^{eff}$ satisfies the relationship $G_c^{eff}{}_{(\xi=2)} < G_c^{eff}{}_{(\xi=1)} < G_c^{eff}{}_{(\xi=0)}$.

For the model in Fig. 4, it is assumed that there is an initial crack inside $\Gamma_0 = [-a_0, +a_0] \times \{0\}$,

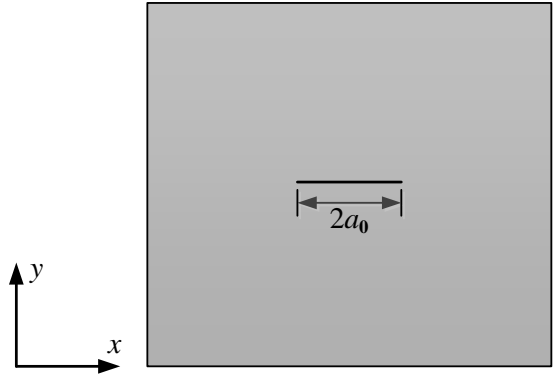

**Fig. 4.** 2D crack ($x = -a_0$ to $x = a_0$)

It can be calculated that the effective initial crack length $a_0^{eff}$ is [47]

$$a_0^{eff} = a_0\left(1+\frac{\pi l}{4a_0\left(h/4c_w l+1\right)}\right) \tag{70}$$

It can be seen from Eq. (70) that when the parameters $h$, $l$ and $G_c$ are the same, $a_0^{eff}$ satisfies the relationship $a_0^{eff}{}_{(\xi=0)} < a_0^{eff}{}_{(\xi=1)} < a_0^{eff}{}_{(\xi=2)}$.

It can be derived from Eq. (69) and (70)

$$a_0^{eff} = a_0\left(1 + \frac{G_c \pi l}{4a_0 G_c^{eff}}\right) \tag{71}$$

It can be seen from the Eq. (71) that $a^{eff}$ and $G_c^{eff}$ are coupled with each other and $G_c^{eff}$ will affect the value of $a^{eff}$. At the same time, it can be found from Table 2 that the stress and strain are also related to the initial crack length and energy release rate. Therefore, when studying the mechanical properties of different damage models, it is necessary to comprehensively consider the influence of the model's effective crack length and effective energy release rate on the fracture properties of materials.

According to the Griffith energy criterion [48,49], the stability condition of the allowable crack increment can be defined

$$\frac{d}{da}(S - \Pi) = 0 \tag{72}$$

Here, $\Pi$ is the additional strain energy caused by the crack,

$$\Pi = \frac{(3-\beta)\pi a_0^2 \left(\sigma_y^{0,b}\right)^2}{8\mu} \tag{73}$$

$S$ is the fracture energy consumed by the crack,

$$S = 2G_c a_0 \tag{74}$$

where $\mu$ is the shear modulus; for the plane strain problem, $\beta = 3 - 4\nu$; for the plane stress problem, $\beta = \frac{3-\nu}{1+\nu}$, and $\sigma_y^{0,b}$ is the tensile stress acting on the solid at infinity.

According to Eqs. (72)-(74), the maximum tensile stress can be derived

$$\sigma_y^{0,b} = \sqrt{\frac{8G_c\mu}{(3-\beta)a_0\pi}} \tag{75}$$

Assuming that the crack will expand only when the stress exceeds the strength of the material and that the crack will open only when the material releases enough energy, the critical stress at failure can be defined as [49]

$$\sigma_y^{0,c} = \sigma_c \sqrt{\frac{2}{1+\sqrt{1+\sigma_c^4\frac{\pi^2 a_0^2(3-\beta)^2}{16G_c^2\mu^2}}}} = \sigma_c \sqrt{\frac{1}{\frac{1}{2}+\sqrt{\frac{1}{4}+\left(\frac{\sigma_c}{\sigma_y^{0,b}}\right)^4}}} \tag{76}$$

Here, to estimate the critical stress load of the material, the peak stress value in Table 2 is taken as $\sigma_c$, and considering the influence of the element size on the stress, the critical stress estimates of different models can be obtained. In the peak stress value $\sigma_c$, $E_0 := E_0$ is the elastic modulus for the plane stress condition and $E_0 := \frac{E_0}{1-v^2}$ for the plane strain condition with Poisson's ratio $v$.

The critical stress scale factor $\eta$ is defined as

$$\eta = \frac{\sigma_c}{\sigma_y^{0,b}} \tag{77}$$

For different models, the critical stress scale factor is

$$\eta = \begin{cases} \dfrac{9}{32}\sqrt{\dfrac{(3-\beta)(1+v)a_0^{eff}\pi}{3l}}, & \xi=0 \\[2mm] \dfrac{1}{4}\sqrt{\dfrac{3(3-\beta)(1+v)a_0^{eff}\pi}{2l}}, & \xi=1 \\[2mm] \sqrt{\dfrac{(3-\beta)(1+v)a_0^{eff}}{2l}}, & \xi=2 \end{cases} \tag{78}$$

As a result, the estimated value of the critical stress load when different models fail is

$$\sigma_y^{0,c} = \sigma_c \sqrt{\dfrac{1}{\dfrac{1}{2} + \sqrt{\dfrac{1}{4} + \eta^4}}} \tag{79}$$

It can be seen from Eq. (79) that the estimated value of the critical stress load when the material fails is not only a function of the material strength and critical energy release rate but also a function of the geometry (initial crack length). Because the theoretical peak stress value is used in the estimation without considering the effect of damage, its value will be too large, which will cause the critical stress proportional coefficient to become larger, and thus the failure stress estimation result decreases. Therefore, the stress in Eq. (79) is a lower limit of the failure stress of the material, and the true value will be greater than this value.

The stress load in the tensile direction is mainly used to estimate the critical stress load, without considering the effect of shear load on the results. Therefore, the analysis results are mainly suitable for the stress estimation of mode I failure and not suitable for the stress analysis of mode II failure.

## 4 Numerical implementation

### *4.1 Finite element discretization*

Because the governing equation of the phase field model in Section 2.2 cannot be expressed in the standard finite element solution format, the user-defined element (UEL) is needed to solve the model. First, the displacement field and its corresponding strain field are discretized

$$u(x) = \sum_{i=1}^{n} N_i u_i = Nu, \quad \varepsilon(x) = \sum_{i=1}^{n} B_i u_i = Bu, \tag{80}$$

where the interpolation matrix is $N := [N_1, \cdots N_i, \cdots, N_n]$ and the geometric matrix is $B := [B_1, \cdots, B_i, \cdots, B_n]$.

Similarly, the phase field and its corresponding gradient field can be discretized into

$$d(x) = \sum_{i=1}^{n} \bar{N}_i d_i = \bar{N}d, \quad \nabla d(x) = \sum_{i=1}^{n} \bar{B}_i d_i = \bar{B}d. \tag{81}$$

The corresponding displacement field and phase field function variation and its derivative are

$$\delta u(x) = \sum_{i=1}^{n} N_i \delta u_i = N\delta u, \quad \delta \varepsilon(x) = \sum_{i=1}^{n} B_i \delta u_i = B\delta u$$

$$\delta d(x) = \sum_{i=1}^{n} \bar{N}_i \delta d_i = \bar{N}\delta d, \quad \nabla \delta d(x) = \sum_{i=1}^{n} \bar{B}_i \delta d_i = \bar{B}\delta d \tag{82}$$

The trial and test spaces are defined as

$$\begin{aligned}
\mathbb{U}_u &:= H_1^E = \{u | u(x) = u^*, \forall x \in \partial \mathcal{B}_u\}, \\
\mathbb{V}_u &:= H_1^0 = \{\delta u | \delta u(x) = 0, \forall x \in \partial \mathcal{B}_u\}, \\
\mathbb{U}_d &:= H_1^E = \{d | d(x) \in [0,1], \dot{d}(x) \geq 0, \forall x \in \Omega\}, \\
\mathbb{V}_d &:= H_1^0 = \{\delta d | \delta d(x) = 0, \forall x \in \Omega\}.
\end{aligned} \tag{83}$$

Then, the following weak form governing equation can be obtained:

$$\begin{cases} \delta u^{\mathrm{T}} \left( \int_{\mathcal{B}} B^{\mathrm{T}} \sigma dV + \int_{\mathcal{B}} b dV \right) = 0 \\ \delta d^{\mathrm{T}} \int_{\Omega} \left[ \bar{N}^{\mathrm{T}} g'(d) \mathcal{H} + \frac{1}{4c_w} \left( \bar{N}^{\mathrm{T}} \frac{1}{l} w'(d) + 2l \bar{B}^{\mathrm{T}} \nabla d \right) \right] dV \geq 0 \end{cases} \tag{84}$$

*4.2 Staggered algorithm implementation*

Due to the strong nonlinearity of the governing equations, it is not easy to converge in the calculation. Here, to improve the convergence of the calculation, the staggered algorithm is used to solve the above control equations of the phase field fracture model [50,51], where in the first iteration at every load step, the history and the phase field are updated for the phase and displacement field elements, respectively, the displacement is based on the phase field value taken from the end of the previous step ($d_n$), and the phase field problem is solved based on the displacement field value in the next step ($\mathcal{H}_{n+1}$). This method is similar to a semi-implicit method, so the convergence of the results is better when solving the governing equations of the phase field model. The element structure of the staggered algorithm in ABAQUS is shown in Fig. 5.

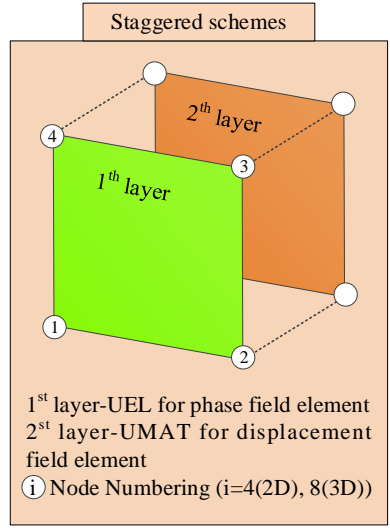

**Fig. 5.** Element structure composition in ABAQUS

Here, the Newton-Raphson iterative method can be used to solve the above equations to obtain the displacement field and the phase field:

$$\begin{bmatrix} K_n^d & 0 \\ 0 & K_n^u \end{bmatrix} \begin{bmatrix} d_{n+1} \\ u_{n+1} \end{bmatrix} = - \begin{bmatrix} r_n^d \\ r_n^u \end{bmatrix} \tag{85}$$

Where each expression is

$$r_n^d = \int_\Omega \left\{ \left[ \frac{G_c}{4c_w l} w'(d) + g'(d)\mathcal{H} \right] (\bar{N})^T + \frac{G_c l}{2c_w} (\bar{B})^T \nabla d \right\} dV$$

$$K_n^d = \frac{\partial r^d}{\partial d} = \int_\mathcal{B} \left\{ \left[ \frac{G_c}{4c_w l} w''(d) + g''(d)\mathcal{H} \right] (\bar{N})^T \bar{N} + \frac{G_c l}{2c_w} (\bar{B})^T \bar{B} \right\} dV \tag{86}$$

$$r_n^u = \int_\mathcal{B} (B)^T \sigma dV - \left( \int_\mathcal{B} (N)^T b dV + \int_{\partial \mathcal{B}} (N)^T t dA \right)$$

$$K_n^u = \frac{\partial r^u}{\partial u} = \int_\Omega (B)^T \left[ \frac{\partial \sigma}{\partial \varepsilon} \right] B dV \tag{87}$$

## 5 Results and discussion

In this section, some examples are used to verify the accuracy of the above theoretical analysis. Due to the similar results of $\xi=1,2$, the examples in this section are used to analyze the mechanical and energy properties of $\xi=0,1$. The first example is used to analyze the linear elastic properties and energy properties of different phase field damage models, and the second example is used to analyze the mechanical properties of these models.

### 5.1 Single element tensile model

The model is a square plate (1×1 mm) with a single quadrilateral element. The Young's modulus of the specimen is set to $E=210\,kN/mm^2$, and the Poisson's ratio is set to $v=0.3$. The critical energy release rate is $G_c=5\times10^{-3}\,kN/mm$, and the internal length scale parameter is $l=0.01mm$. The applied tensile load is $u=0.1mm$. The model and boundary conditions are shown in Fig. 6, and the applied load method is shown in Fig. 7.

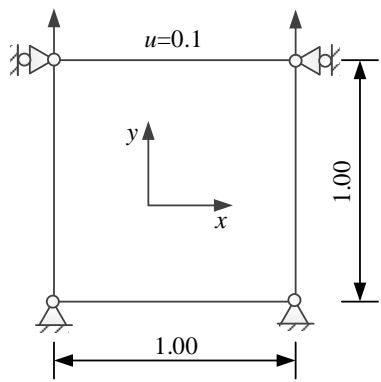

**Fig. 6.** Uniaxial tensile model

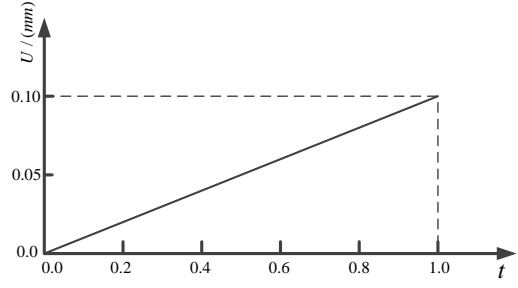

**Fig. 7.** Loading mode

To simplify the calculation, it is assumed that:

$$\begin{cases} \varepsilon_{22} \neq 0 \\ \varepsilon_{11} = \varepsilon_{12} = 0 \end{cases} \quad (88)$$

The stresses can be calculated directly: $\sigma_{22}^0 = E_{22}\varepsilon_{22}$, where $E_{22}$ is the element of the plane strain stiffness matrix: $E_{22} = \dfrac{E(1-\nu)}{(1+\nu)(1-2\nu)}$, then the total elastic energy density is

$$\psi_0 = \frac{1}{2} E_{22} \varepsilon_{22}^2 \quad (89)$$

The degradation function is taken as $g(d) = (1-d)^2$; in case $\xi = 0$, the critical strain energy, damage variable and critical strain can be calculated as

$$\begin{aligned} &\mathcal{H}_c = 0 \\ &d = \frac{2l\psi_0}{G_c + 2l\psi_0} = \frac{E_{22}\varepsilon_{22}^2 l}{G_c + E_{22}\varepsilon_{22}^2 l} \\ &\varepsilon_c = 0.02428 \end{aligned} \quad (90)$$

For model $\xi = 0$ with additional critical strain energy $\mathcal{H}_c$, the critical strain energy, damage variable and critical strain can be calculated as

$$\mathcal{H}_c = \frac{G_c}{2l}$$

$$d = \frac{\langle E_{22}\varepsilon_{22}^2 l - G_c \rangle_+}{G_c + \langle E_{22}\varepsilon_{22}^2 l - G_c \rangle_+} \tag{91}$$

$$\varepsilon_c = 0.04206$$

Similarly, in case $\xi = 1$, the critical strain energy, damage variable and critical strain can be calculated as

$$\mathcal{H}_c = \frac{3G_c}{16l}$$

$$d = \frac{\langle 8E_{22}\varepsilon_{22}^2 l - 3G_c \rangle_+}{8E_{22}\varepsilon_{22}^2 l} \tag{92}$$

$$\varepsilon_c = 0.02575$$

Then, the stress along the $y$ direction can be calculated as $\sigma_{22} = (1-d)^2 \sigma_{22}^0$.

According to the loading method shown in Fig. 7, the simulation results of the damage field $d$ and stress $\sigma_{22}$ can be obtained according to the staggered algorithm and compared with the analytical solution. The results are shown in Figs. 8 and 9.

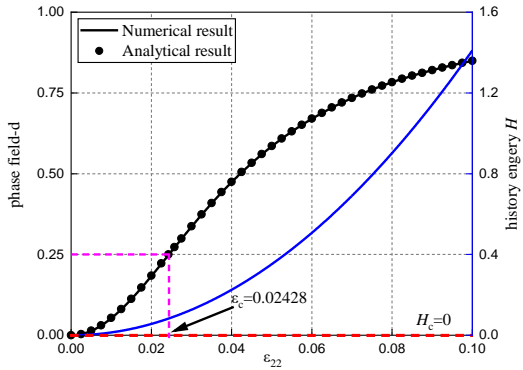 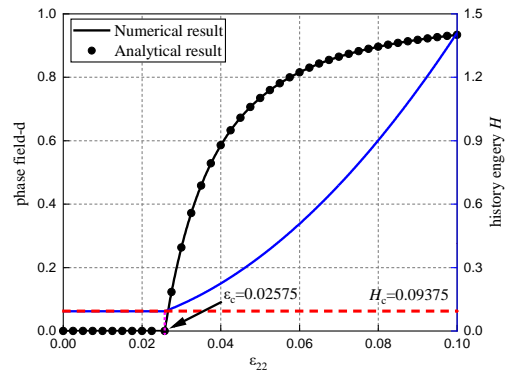

(a) $\xi = 0$ (without critical strain energy)   (b) $\xi = 1$

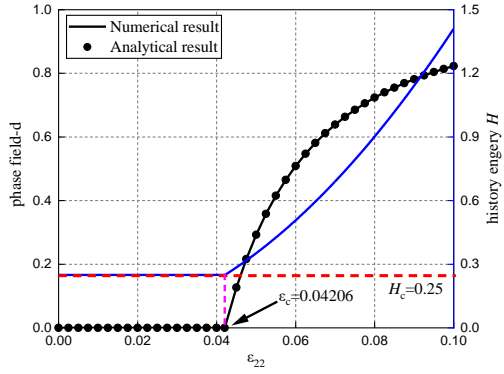

(c) $\xi = 0$ (with critical strain energy)

**Fig. 8.** Damage field $d$; the black solid line is the numerical simulation result of the damage field, the black circle dot is the analytical solution of the damage field, the blue solid line is the historical strain energy, and the red dashed line is the critical strain energy.

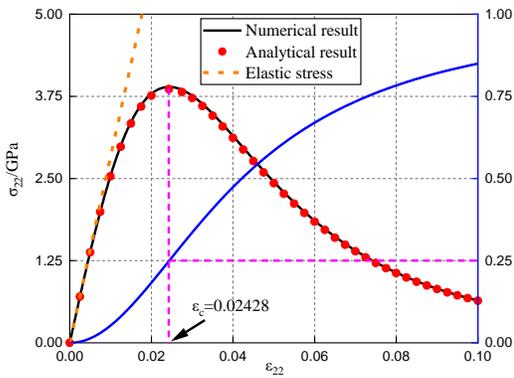

(a) $\xi = 0$ (without critical strain energy)

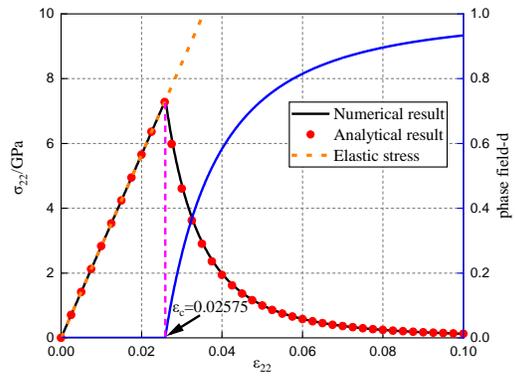

(b) $\xi = 1$

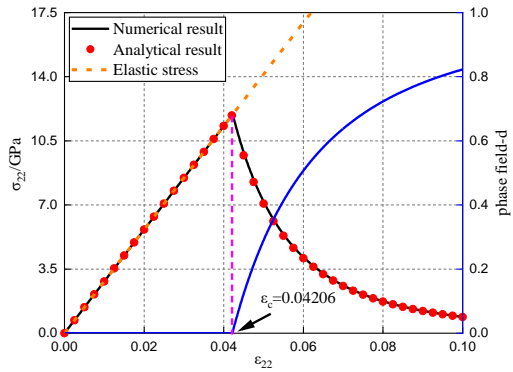

(c) $\xi = 0$ (with critical strain energy)

**Fig. 9.** Stress $\sigma_{22}$; the black solid line is the numerical simulation result of the stress, the red circle dot is the analytical solution of the stress, the blue solid line is the damage field, and the orange dashed line is the numerical simulation result of elastic stress $\sigma_{22}^0$.

As shown in Fig. 8(a), in terms of model $\xi = 0$, when the load starts to be applied, damage begins to occur, and the historical strain energy also increases from the critical strain energy $\mathcal{H}_c = 0$. The damage field increases to the critical value $d = 0.25$ as the strain increases to the critical strain. However, for model $\xi = 1$, as shown in Fig. 8(b), when the load starts to load, the damage does not start immediately and keeps $d = 0$, and the historical strain energy of this model remains at the critical value $\mathcal{H}_c = 0.09375$. With the continuous increase in the load, when the strain increases to the critical strain, the strain energy begins to increase, and damage begins to occur. By comparison (Fig. 8(a) and (b)), it can be found that the damage field begins to occur when the model is loaded because there is no control of the critical strain energy for model $\xi = 0$, as shown in Eq. (62); however, damage does not occur for model $\xi = 1$ when the strain energy of the model is less than or equal to the critical strain energy. Due to the existence of the critical strain energy $\mathcal{H}_c$, the damage always remains $d = 0$. Only when the strain energy is greater than the critical strain energy does damage begin to occur, and the results are consistent with Eq. (62). However, when a critical strain energy $\mathcal{H}_c$ is introduced into the model $\xi = 0$ (Fig. 8(c)), damage does not occur immediately after the load starts to load, but when the strain energy is greater than the critical strain energy, damage begins to occur. This analysis verifies that the critical strain energy affects the starting time of the damage.

Fig. 9(a) shows that when the load starts to load, damage occurs immediately for model

$\xi = 0$, leading to a decrease in the stress increment. When the strain reaches the critical strain, the damage field increases to the critical value $d=0.25$, and the stress increases to the limit value; then, the stress decreases with increasing load. However, for model $\xi = 1$, as shown in Fig. 9(b), when the load starts to be applied, the stress begins to increase, but its increment does not change, so that the stress-strain maintains a linear elastic change. Only when the strain increases to the critical strain does damage start to occur. At this time, the stress begins to decrease with increasing load. By comparing models $\xi = 0$ and $\xi = 1$ (Fig. 9(a) and (b)), it can be found that the stress-strain relationship obtained by model $\xi = 0$ lacks the linear elastic stage and does not meet the linear elastic fracture mechanics results, while the results obtained by model $\xi = 1$ meet the linear elastic relationship at the beginning of the material change stage, which is consistent with the linear elastic fracture mechanics, indicating that model $\xi = 1$ is more accurate than model $\xi = 0$ in simulating the brittle failure process of linear elastic materials. However, when the critical strain energy is introduced into the model $\xi = 0$ (Fig. 9(c)), the obtained stress-strain relationship satisfies the linear elastic relationship at the initial stage and has the elastic limit stress. Through this example, it is further verified that the critical strain energy $\mathcal{H}_c$ affects the time when damage occurs and further affects the stress-strain change process so that the model has the linear elastic stage in the initial change process and has the elastic limit stress, which shows the accuracy of the theoretical analysis in Section 3.2.

*5.2 Single edge notched tensile model*

In this section, the model is a square plate (1×1 mm), there is an initial crack $\Gamma$ in the middle area of the model, with a length of $a_0 = 0.5mm$, the lower boundary is fixed, and the tensile load is applied at the upper boundary. The boundary condition and load form of the model are shown in Fig. 10, and the material parameters are elastic modulus $E = 210\,kN/mm^2$ and Poisson's ratio $v=0.3$. Two different cases are given to study the effects of different parameters on the mechanical properties of materials. The first case is to study the mechanical properties of different phase field models when the effective critical energy release rate $G_c^{eff}$ changes under the same effective initial crack length $a_0^{eff}$; the second case is to study the mechanical properties of different phase field models when the same effective initial crack length $a_0^{eff}$ changes under the same conditions of effective critical energy release rate $G_c^{eff}$. In both cases, the model element size is the same $h = 0.005\,mm$.

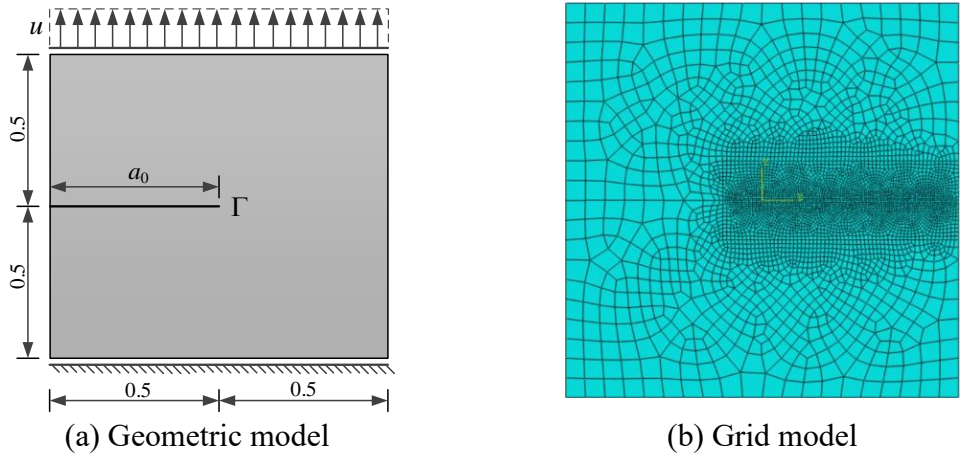

(a) Geometric model  (b) Grid model

**Fig. 10.** Single edge notched model

*5.2.1 Case I*

For the first case, the model parameters are shown in Table 3. According to Eq. (70), the effective initial crack length is calculated as $a_0^{eff}{}_{(\xi=0)} = a_0^{eff}{}_{(\xi=1)} = 0.50628\,mm$, and according to Eq. (69), the effective critical energy release rate is calculated as $G_c^{eff}{}_{(\xi=0)} = 0.003375\,kN/mm$, $G_c^{eff}{}_{(\xi=1)} = 0.003230\,kN/mm$, and $G_c^{eff}{}_{(\xi=1)} < G_c^{eff}{}_{(\xi=0)}$. The mechanical properties of different phase field models are studied when $G_c^{eff}$ changes. The simulation results are shown in Fig. 11.

| Parameters | $\xi = 0$ | $\xi = 1$ |
|---|---|---|
| $G_c$/(kN/mm) | $2.7\times10^{-3}$ | $2.7\times10^{-3}$ |
| $l$/mm | 0.01 | 0.00956 |

**Table 3.** Model parameter

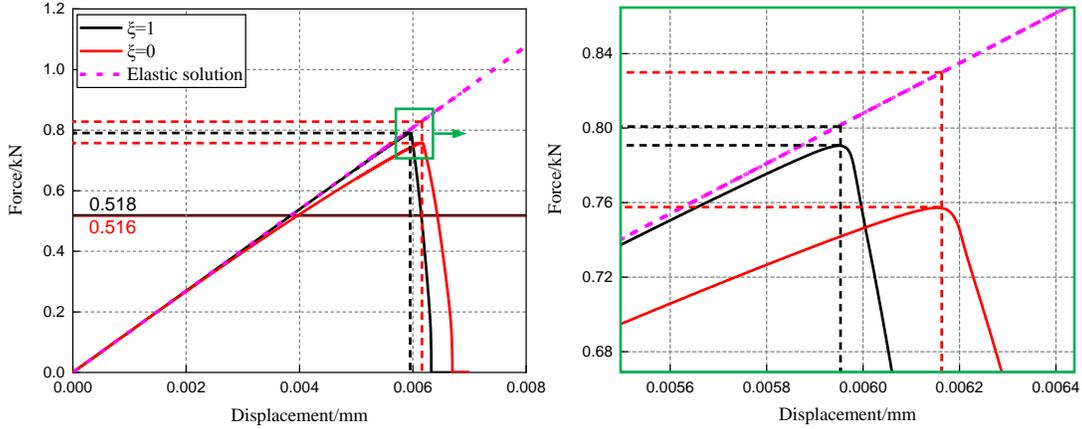

**Fig. 11.** Displacement-load curve

Fig. 11 shows that the crack opening displacement of model $\xi = 0$ is greater than that of model $\xi = 1$, but the material strength at failure is less than that of model $\xi = 1$. This is because the effective initial crack length $a_0^{eff}$ of different models is the same in this example, but the effective critical energy release rate of $G_c^{eff}$ of model $\xi = 0$ is greater than that of model $\xi = 1$, which improves the fracture toughness and delays the destruction process of the material. Therefore, the crack opening displacement of model $\xi = 0$ is larger than that of model $\xi = 1$.

In the process of material failure, because 25% of the damage has occurred in the material when the stress of model $\xi = 0$ reaches the maximum stress, the material strength at failure of model $\xi = 0$ is less than the failure strength of model $\xi = 1$. If the material in model $\xi = 0$ is not damaged at the maximum stress, the material failure strength of model $\xi = 0$ will be greater than that of model $\xi = 1$ under the condition that the effective critical energy release rate $G_c^{eff}$ of model $\xi = 0$ is greater than that of model $\xi = 1$, as shown in the red extended dashed line in Fig. 11.

The critical stress value of material failure can be estimated by Eq. (79), and then the force load of material failure can be calculated, as shown in Fig. 11. The estimated results of the two models are very similar, but they are less than the simulation results, which is consistent with the theoretical analysis results in Section 3.3. The estimated stress given in Eq. (79) is a lower limit of the material failure stress, and the actual failure strength is greater than the estimated value, which verifies the accuracy of the theoretical analysis in Section 3.3.

*5.2.2 Case II*

For the second case, the model parameters are shown in Table 4. According to Eq. (69), the effective critical energy release rate is calculated as $G_{c\ (\xi=0)}^{eff}=G_{c\ (\xi=1)}^{eff}=0.0034375\ kN/mm$, and according to Eq. (70), the effective initial crack length is calculated as $a_{0\ (\xi=0)}^{eff}=0.50628mm$, $a_{0\ (\xi=1)}^{eff}=0.50286mm$, and $a_{(\xi=1)}^{eff}<a_{(\xi=0)}^{eff}$. The mechanical properties of different phase field models are studied when $a_0^{eff}$ changes. The simulation results are shown in Fig. 12.

| Parameters | $\xi=0$ | $\xi=1$ |
|---|---|---|
| $G_c$/(kN/mm) | 2.75×10⁻³ | 2.5×10⁻³ |
| $l$/mm | 0.01 | 0.005 |

**Table 4.** Model parameter

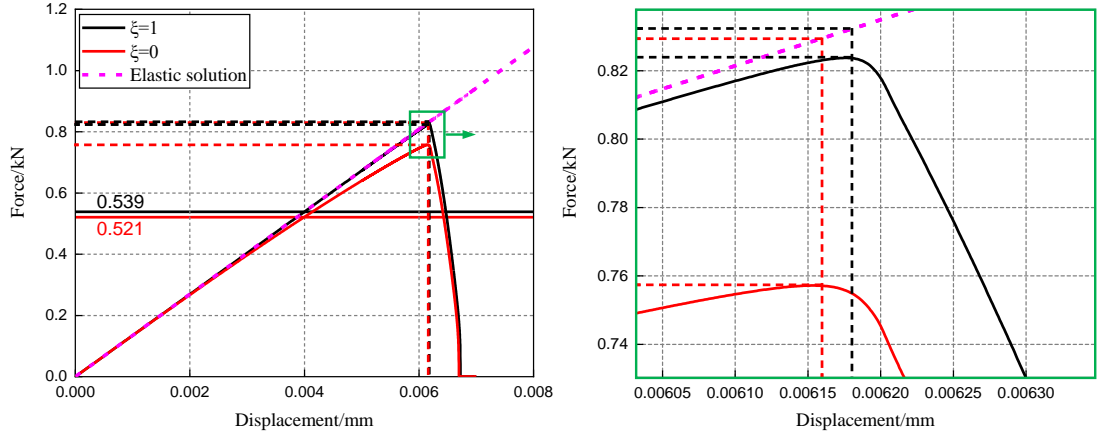

**Fig. 12.** Displacement-load curve

It can be seen from Fig. 12 that the crack opening displacement of the model $\xi = 0$ is approximately equal to the crack opening displacement of the model $\xi = 1$, and the material strength at failure is less than the failure strength of the model $\xi = 1$. The reason is that the effective critical energy release rate $G_c^{eff}$ of different models is the same in this example. Although the effective initial crack length $a_0^{eff}$ of the model $\xi = 0$ is greater than that of the model $\xi = 1$, making its effective load-bearing length smaller, the crack opening displacement of the model $\xi = 0$ is approximately equal to that of the model $\xi = 1$ under the combined influence of the two factors. The results show that the effect of the effective critical energy release rate on the stress is greater than that of the effective initial crack length.

In the process of material failure, because 25% of the damage has occurred in the material when the stress of model $\xi = 0$ reaches the maximum stress, the material strength at failure of model $\xi = 0$ is less than the failure strength of model $\xi = 1$. If the material in model $\xi = 0$ is not damaged at the maximum stress, the material failure strength of model $\xi = 0$ will be approximately equal to that of model $\xi = 1$ under the effective critical energy release rate

$G_c^{eff}$ of different models is the same, as shown in the red extended dashed line in Fig. 12.

The critical stress of the material at failure can be estimated by Eq. (79), and then the force load at failure can be calculated as shown in Fig. 12. The estimated results of the two models are similar, but they are less than the simulation results, which is similar to the results and reason of the first case and further verifies the accuracy of the theoretical analysis in Section 3.3.

## 6 Conclusion

In this paper, the energy and mechanical properties of various phase field damage models are analyzed, and the accuracy of the theoretical analysis is verified by some examples. Based on our results, the most important findings are as follows:

(1)  The properties of the homogenous solution and nonhomogenous solution of the governing equation of the phase field model, as well as the stress function, are studied. Under special boundary conditions, the nonhomogenous solution is the geometric crack function of the phase field model.

(2)  The critical strain energy of the phase field model affects the damage, stress and strain properties of the material. When there is a critical strain energy in the phase field model, damage does not occur in the initial stage of loading; only when the elastic strain energy increases to the critical strain energy does damage begin to occur. Therefore, the stress-strain relationship of the model has a linear elastic stage, and the elastic limit stress exists in the model. However, when there is no critical strain energy in the damage model, the model will be damaged at the beginning of loading, resulting in a nonlinear change in the stress-strain relationship from the beginning of loading, and there is no elastic limit stress in the model.

(3) In the process of material failure, the effective critical energy release rate has a greater impact on the failure stress of the material than that of the effective initial crack length. When the effective critical energy release rate or the effective initial crack length is the same, the failure stress of model $\xi=0$ (AT2) is less than that of model $\xi=0$ (AT1) due to 25% of the damage occurring at the maximum stress of model $\xi=0$ (AT2).

(4) An estimated value of the material failure stress of the phase field damage model is given in this paper, which is a lower limit of material failure stress, and the actual failure strength will be greater than the estimated value.

**Acknowledgements**

Financial support from National Key Program of Basic Research on Foundation Strengthening Project（No.2020-JCJQ-ZD-076-02）.

**Declaration of competing interest**

The authors declare that they have no known competing financial interests or personal relationships that could have appeared to influence the work reported in this paper.